\input amstex
\documentstyle{amsppt}
\magnification 1200
\NoBlackBoxes
\hsize= 14.6 truecm
\vsize= 19 truecm
\hoffset= .5 truecm
\voffset= 2.1 truecm

\def\ns{\hskip -.15 truecm}

\def\main{Theorem 1 }

\def\notation{Notation 2.1 }
 \def\smallemma{Lemma 2.2 }
  \def\superlemma{Lemma 2.3 }
 \def\oddcurve{Proposition 2.4 }
\def\oddformulaone{(2.5.1) }
\def\oddformulatwo{(2.5.2) }
\def\oddmain{Theorem 2.6 }
 \def\oddremark{Remark 2.7 }
 \def\defNo{Definition 2.8 }
 \def\oddNo{Theorem 2.9 }
 \def\evencurve{Proposition 2.10 }
 \def\evenmain{Theorem 2.11 }
 \def\evenremark{Remark 2.12 }
 \def\evenNo{Theorem 2.13 }
 \def\evenNoquestion{Question 2.14 }

 \def\genstruct{Theorem 3.1 }
  \def\scroll{Theorem 3.2 }
  \def\flatstruct{Theorem 3.3 }
  \def\flatone{(3.3.1) }
  \def\flattwo{(3.3.2) }
  \def\examples{Proposition 3.4 }
  \def\Greenremark{Remark 3.5 }

\def\Catone{[Ca1] }
\def\Cattwo{[Ca2] }
\def\Ci{[Ci] }
\def\Green{[G] }
   \def\Miranda{[M] }
    \def\HM{[HM] }
    \def\GP{[GP] }

\topmatter
\title Triple canonical covers of varieties of minimal degree
\endtitle
\author Francisco Javier Gallego \\ and \\ Bangere  P. Purnaprajna
\endauthor
\abstract{In this article we study pluriregular varieties $X$ of general
  type with base-point-free canonical bundle whose canonical morphism
  has degree $3$ and maps $X$ onto a variety of minimal degree $Y$. We
  carry out our study from two different perspectives. 

First we study in Section 2
  the canonical ring of $X$ describing completely the degrees of its
  minimal generators. We apply this to the study of the projective
  normality of the images of the pluricanonical morphisms of $X$. 
 Our study of the canonical ring of $X$ also shows that, if the
  dimension of $X$ is greater than or equal to $3$, there does not exist a
  converse to a  theorem of M. Green that bounds the degree of the 
generators of the canonical ring of $X$. 
This is in sharp contrast with the situation in dimension $2$ where such
converse exists, as proved by the authors in a previous work.

Second, we study in Section 3, the structure of the canonical morphism
of $X$. We use this to show among other things the non-existence of some a
priori plausible examples of triple canonical covers of varieties of
minimal degree. We also characterize the targets of flat canonical covers
of varieties of minimal degree. Some of the results of Section 3 are 
more general and apply to varieties $X$ which are not necessarily
regular, and to targets $Y$ that are scrolls which are not of minimal
degree.    
}
\endabstract
\address{Francisco Javier Gallego: Dpto. de \'Algebra,
 Facultad de Matem\'aticas,
 Universidad Complutense de Madrid, 28040 Madrid,
Spain}\endaddress
\email{FJavier\_Gallego\@mat.ucm.es}\endemail
\address{ B.P.Purnaprajna:
405 Snow Hall,
  Dept. of Mathematics,
  University of Kansas,
  Lawrence, Kansas 66045-2142}\endaddress
\email{purna\@math.ukans.edu}\endemail
\endthanks
\thanks{The first author was partially supported by MCT project
number BFM2000-0621 and by UCM project number PR52/00-8862.
The second author was partially supported by the
General Research Fund of the University of Kansas at
Lawrence. The first author is grateful for the hospitality of the
Department of Mathematics of the University of Kansas at
Lawrence.}\endthanks

\endtopmatter
\document

\vskip .3 cm

\headline={\ifodd\pageno\rightheadline \else\leftheadline\fi}
\def\rightheadline{\tenrm\hfil \eightpoint TRIPLE CANONICAL COVERS OF
  VARIETIES OF MINIMAL DEGREE
 \hfil\folio}
\def\leftheadline{\tenrm\folio\hfil \eightpoint F.J. GALLEGO \&
  B.P. PURNAPRAJNA \hfil}

\heading 1. Introduction \endheading

The purpose of the article is to study pluriregular varieties (these are,
according to Catanese, varieties such that the intermediate cohomology of
its structure sheaf vanish) of general type $X$ with base-point-free
canonical divisor and such that the canonical morphism $\pi$ is a
generically $3:1$ morphism onto a variety of minimal degree $Y$.  The
varieties of general type whose canonical morphism maps to a variety of
minimal degree play an important role in various settings.  In dimension
$2$ they are central in the classification of surfaces of general type
with small  $c_1^2$ and in questions about degenerations and the moduli
of surfaces of general type. This is illustrated by the works of Horikawa,
Catanese, Konno, Mendes Lopes and Pardini, among others. In higher
dimensions, pluriregular varieties of general type have been studied by 
F. Catanese (see \Catone and \text{\Cattwo \ns)} and M. Green (\Green \ns) 
who prove results on the canonical ring 
of these varieties.

\bigskip

We study these varieties from two different perspectives. The first
perspective we take is the study of the structure of their canonical
ring. This is handled in Section 2. More precisely our main goal is 
to find out the degrees of the generators. This goal is accomplished 
in \oddmain for odd dimensional varieties and \evenmain for even 
dimensional varieties, where we obtain a full description of the degrees 
in terms of the dimension of $X$ and the degree of $Y$. We put these 
two results together in the following

\proclaim{\main} 
Let $X$ be a smooth, pluriregular variety of general type
of dimension $m \geq 2$ 
and let $Y 
\subset
\bold P^{r+m-1}$ be an $m$-fold of minimal degree. Assume that 
$K_X$ is base-point-free and let
$\pi: X \longrightarrow Y$ be the canonical morphism of $X$. 
Assume that $\pi$ is a generically finite morphism of degree $3$. Then

\item{1)} if $m$ is odd,  
the canonical ring of $X$ is generated by its part of
degree
$1$ and
$\frac{r}{2}$ linearly indendent generators in degree $\frac{m+1}{2}$;

\item{2)} if $m$ is even, 
the canonical ring of $X$ is generated by 
its elements of degree $1$, $r$ linearly independent  
generators in  degree
$\frac{m+2}{2}$ and $r-1$ linearly independent generators in degree $m+1$. 

\endproclaim

Apart from its intrinsic interest, \main  has several applications.  In
the first place the degrees of the generators of the canonical ring are
related to the degrees of the generators of the pluricanonical rings and
therefore, with the very ampleness and projective normality associated to
pluricanonical morphisms on the variety $X$. In this regard, using 
\main (1) for the odd dimensional case, we give a necessary and 
sufficient condition for multiples of the canonical bundle $K_X$ on 
$X$ to give projectively normal embeddings. Precisely (see \oddNo \ns)
we obtain

\bigskip

{\it Let $X$ an odd dimensional variety satisfying
the hypothesis of \main \ns. Assume in addition that $K_X$ is ample. 
Then the line
bundle
$K_X^{\otimes n}$ is very ample and the image of the morphism it induces
is projectively normal if and only if $n \geq \frac{m+1}{2}$.} 

\bigskip

In the even dimensional case using \main (2) we obtain a
sufficient condition for a multiple of $K_X$ to induce a projectively
normal embedding (see \evenNo \ns). 

\bigskip

Another interesting consequence of \main concerns this general 
question: ``what is the minimal degree of the generators of the 
canonical ring of a variety of general type?.''  Green (see \Green \ns) 
proved that the canonical ring of a pluriregular variety of general 
type is generated in degree less than or equal to $n$ if the canonical 
morphism does not map $X$ onto a variety of minimal degree. 
Here is yet another instance where varieties of minimal degree 
play an important role. For surfaces, the above mentioned result 
by Green was also obtained independently by Ciliberto (see \Ci \ns). 
In \GP we proved
a converse of this result of Green and Ciliberto. Precisely, 
we found out the degrees of the generators of the canonical ring of a 
regular surface of general type $X$ whose canonical morphism mapped $X$ 
to a surface of minimal degree. As a corollary of this we proved 
(cf. \GP \ns, Corollary 2.8) a converse to the result of Ciliberto and Green.
For higher dimensions, \main  and the examples at the end of Section 3 
draw a very different picture. Indeed \main  and \examples show that a
converse to Green's result, that is, a result saying that the canonical 
ring of a pluriregular variety of general type is generated in degree less 
than or equal to $n$ {\it if and only if} 
the canonical morphism does not map $X$ onto a variety of
minimal degree, is false in dimension greater than $2$. The reason the 
converse fails is due to the fact that for the varieties studied in 
Section 2, the canonical ring is generated in degree much lower than the 
Green's bound, roughly, in degree half of Green's bound. This fact is noted 
in \Greenremark \ns.

\bigskip

The second approach we take is to study the structure of the canonical
morphism of varieties $X$ of general type which are triple covers of  
varieties of minimal degree $Y$. This is dealt with in Section 3. 
We will mention here some of the interesting consequences of this
study. 
We obtain results which show that the fact of working with a morphism 
induced by a canonical subseries imposes  strong restrictions on the parity 
of the dimension of $X$ and $Y$ and of the degree of $Y$. These results hold 
for varieties of general type which are more general than those considered 
in Section 1. For instance, as a consequence of \text{\scroll \ns,} we obtain 
this result which holds for canonical morphisms of arbitrary degree and 
for any scroll $Y$ not necessarily fibered over
$\bold P^1$: 

\medskip

(1) An even dimensional variety of general type does not
admit a generically finite canonical morphism of odd degree to a
scroll. 

\medskip

When the canonical morphism is generically finite of degree $3$, 
we go further in our study and obtain \genstruct and \flatstruct \ns. 
 We mention 
below some of the consequences of these theorems:

\medskip

(2) We show that there is a relation between the dimension   
$m$ of $X$ and the degree $r$ of $Y$. More precisely, we see that the
condition of $m$ being odd forces $r$ to be even. In particular there are
no odd dimensional varieties of general type which are triple canonical
covers of linear $\bold P^m$.

\medskip 

With the additional hypothesis of flatness on the canonical morphism, we
obtain the following stronger results: 

\medskip

(3) A variety of general type which is a canonical triple cover of a
variety of minimal degree is pluriregular (see \genstruct \ns). 

\medskip

(4) \flatstruct shows that all triple flat canonical
covers $X @>\pi >> Y$  of varieties of minimal degree are 
such that $\pi_*\Cal O_X$ splits as
vector bundle over $\Cal O_Y$ in the same way as a cyclic cover.

\medskip

(5) We determine precisely what are the varieties of
general type which occur as targets of canonical morphisms of
degree $3$. This is achieved by \flatstruct and \examples \ns. 

\medskip

Moreover, the pluriregular varieties of general type 
constructed in \examples illustrate the theorems that 
we prove in this paper.

\bigskip

We will say now a few words on the techniques used to prove the results of
Section 2. These techniques involve the study of the $\Cal {O}_{Y}$-algebra
structure on $\pi _{*}\Cal{O}_{X}$ to find the multiplicative
structure of the canonical ring of $X$. 
Even though we reduce the problem
from a complicated variety to a simpler variety, as a variety of minimal
degree $Y$ appear to be when compared to the variety $X$, 
there are natural difficulties that arise in the
process. The proof of \oddmain and \evenmain 
involves the study of multiplication of global
sections of line bundles on the $m$-dimensional 
variety of general type $X$. To do so we
relate the problem of multiplication of global sections of line bundles on
$X$ to that of multiplication of global
sections of line bundles on a curve $C$ which is an $m-1$ 
complete intersection on $X$. We achieve this 
through an inductive process that we develop in \superlemma \ns. An impression
that things become simple on curves is misleading since by a process such
as this one obtains bundles of small degree on $C$ 
which are difficult to handle in the context of multiplication maps. 
To circumvent this we transfer the
problem of multiplication maps on $C$ to a rational normal curve 
obtained as $m-1$ complete intersection on $Y$ by means of pushing down
from $C$ by $\pi$.  
We have traded one difficulty, that of handling small degree
bundles on the curve $C$ of $X$, to another, namely, the study 
of certain maps on $\bold P^1$. These maps involve 
global sections of vector bundles on $\bold P^1$ but 
are not multiplication maps of global sections. 
The key ingredient to this  approach, 
which allows us to
settle the problem, is to use the
algebra structure of $\Cal
O_C$ over $\Cal O_{\bold P^1}$ to interpret the 
maps appearing on $\bold P^1$
and to see their relation with the multiplication maps on $C$.
Note finally that, unlike double covers, a triple cover is not 
determined by the branch locus and the algebra structure it induces 
on $\Cal O_X$ might be quite complicated. This is another difficulty one
encounters but we are able to overcome it in the context of this paper. 

\medskip

{\bf Convention.} Throughout this article the base field $\bold k$
has characteristic $0$ and $\overline{\bold k}=\bold k$. By variety we
will mean an irreducible projective variety over $\bold k$.

       \heading 2. The canonical ring of triple covers of varieties of
       minimal degree
       \endheading

        The purpose of this section is to study the canonical ring of
        a  pluriregular  variety
        $X$ of general type with base-point-free canonical bundle $K_X$ and
        such that its canonical morphism maps $X$ generically $3:1$ onto a
        variety of minimal degree.  Precisely we want to prove \main of the
        introduction, which tells the
        degrees of the generators of the canonical ring of $X$. 
        To do so
        we will need to establish some notation and to prove some auxiliary
        lemmas and results. 

        \bigskip

         {\bf \notation \ns.}
          Let $X$ be a smooth variety of general
          type of dimension $m$ and let $K_X$ be its canonical bundle.
          Assume that $K_X$ is base-point-free. Let $X_1 \subset \cdots
          \subset X_{m'} \subset \cdots X_m=X$ be smooth irreducible
          subsequent $m'$-dimensional, complete intersections of members
          of $|K_X|$. We will also denote $X_1$ by $C$. We denote
          $K_X|_{X_{m'}}$ by
          $L_{m'}$ (in particular $L_m=K_X$) and we will also denote $L_1$
          by $\theta$. Finally let
          $$H^0(X_{m'},L_{m'}^{\otimes s}) \otimes
          H^0(X_{m'},L_{m'}^{\otimes t}) @> \alpha(s,t;m') >>
         H^0(X_{m'},L_{m'}^{\otimes s+t}) $$
           be the usual multiplication map of global sections. 

           We will use also some abridged notation:

            \item{} $\alpha(s,t;m)=\alpha(s,t)$,
            \item{} $\alpha(s,t;1)=\beta(s,t)$,
            \item{} $\alpha(s,1;m')=\alpha(s;m')$
            \item{} $\alpha(s,1;m)=\alpha_s$ and
            \item{} $\alpha(s,1;1)=\beta_s$.

               \bigskip

               The idea of the proofs of \main is to relate the
                computation of generators of the canonical ring to the knowledge of
                the images of certain maps of multiplication of global sections of
                line bundles on $X$. We will find out the image of these maps by
                studying analogous maps of multiplication of global sections of line
                bundles on the varieties $X_{m'}$ defined above, and eventually, by
studying maps on the curve $C$. In order to
realize the link between the maps on $X_{m'}$ and the maps on $X_m$ we
need a couple of lemmas.

 \proclaim{\smallemma} Let $X$ be a smooth,  pluriregular 
 variety of general
 type of dimension $m$, such that its canonical bundle $K_X$  is
 base-point-free. 
 Then:

 \item{1)} $H^b(L_{m'}^{\otimes a})=0$ for all $2 \leq m' \leq m$,
  $1
  \leq b
  \leq m'-1$
  
  \item{2)} The natural maps of restriction of global sections
   $$H^0(L_{m'}^{\otimes n}) \longrightarrow H^0(L_{m'-i}^{\otimes
   n})$$ surject for all  $2 \leq m' \leq m$, all $1 \leq i \leq m'-1$ and
   all $n \geq 1$. In particular, the maps
   $$\displaylines{H^0(L_{m'}^{\otimes n}) \longrightarrow
   H^0(L_{m'-1}^{\otimes n}) \text{ and }
   \cr H^0(K_X^{\otimes n}) \longrightarrow
   H^0(\theta^{\otimes n})}$$ surject for all $n \geq 1$.

   \endproclaim

    \noindent{\it Proof.} If $m'=m$, 1) follows from
     hypothesis, Kawamata-Viehweg vanishing and Serre duality. Arguing by
     induction on the codimension and using the long exact sequence
     of cohomology which arises from
     $$0 \longrightarrow  L_{m'}^{\otimes n-1}
     \longrightarrow L_{m'}^{\otimes n} \longrightarrow
     L_{m'-1}^{\otimes n} \longrightarrow 0$$ 
     we prove 1) for all $2 \leq m' \leq m$.
     To prove 2) it is enough to show the surjectivity of
     $$H^0(L_{m'}^{\otimes n}) @> \mu >> H^0(L_{m'-1}^{\otimes n})$$
     for all $2 \leq m' \leq m$ and $n \geq 1$, since the maps
     $$H^0(L_{m'}^{\otimes n})
     \longrightarrow H^0(L_{m'-i}^{\otimes n})$$ are composite of maps
     like the previous one. But the surjectivity of $\mu$ follows at
     once from the vanishing of
     $H^1(L_{m'}^{\otimes n-1})$, which follows from 1). \qed

        \proclaim{\superlemma}   Let $X$ be a smooth, pluriregular
        variety of general
        type of dimension $m$, such that its canonical bundle $K_X$  is
        base-point-free. 

        \item{1)} Let $s_1, s_2 \geq 1$. If $$ \displaylines{ H^0(K_X^{\otimes
               s_1}) \otimes
               H^0(K_X^{\otimes s_2}) @> \alpha(s_1,s_2;m) >>
               H^0(K_X^{\otimes s_1+s_2}) \text{ surjects, so does }\cr
               H^0(\theta^{\otimes s_1}) \otimes H^0(\theta^{\otimes s_2})
               @> \alpha(s_1,s_2;1) >> H^0(\theta^{\otimes s_1+s_2}) \ .}$$

               \item{2)} Given $s_1, s_2 \geq 1$, assume that the maps
                $$ H^0(\theta^{\otimes s'})
                \otimes H^0(\theta^{\otimes s_2})
                @> \alpha(s',s_2;1) >> H^0(\theta^{\otimes s'+s_2})$$ surject for
                all $1
                \leq s' \leq s_1$. Then the maps  $$ H^0(L_{m'}^{\otimes s'})
                \otimes H^0(L_{m'}^{\otimes s_2})
                @> \alpha(s',s_2; m') >> H^0(L_{m'}^{\otimes s'+s_2})
                $$ 
                 surject for all
                $1
                \leq s' \leq s_1$ and all $1 \leq m' \leq m$.

                \item{3)} Given $s_1, s_2 \geq 1$, assume that the maps
                 $$ H^0(\theta^{\otimes s'})
                 \otimes H^0(\theta^{\otimes s_2})
                 @> \alpha(s',s_2;1) >> H^0(\theta^{\otimes s'+s_2})$$ 
surject for
                 all $1
                 \leq s' \leq s_1-1$. Then the image of the map 
 $$H^0(L_{m'}^{\otimes s_1})
                 \otimes H^0(L_{m'}^{\otimes s_2})
                 @> \alpha(s_1,s_2;m') >> H^0(L_{m'}^{\otimes s_1+s_2})$$ is the inverse
                 image of the image of map
                 $$H^0(\theta^{\otimes s_1})
                 \otimes H^0(\theta^{\otimes s_2})
                 @> \alpha(s_1,s_2;1) >> H^0(\theta^{\otimes s_1+s_2})$$
                 by the obvious restriction map
                 $$H^0(L_{m'}^{\otimes s_1+s_2}) \longrightarrow H^0(\theta^{\otimes
                   s_1+s_2}) \ . $$
                   Moreover, the codimension of $\alpha(s_1,s_2;m')$ in
                   $H^0(L_{m'}^{\otimes s_1+s_2})$ is the same for all $1 \leq m' \leq
                   m$.

                   \endproclaim

                    \noindent{\it Proof.} First we prove 1).
                     Consider the following commutative diagram:

                     $$\matrix
                       H^0(K_X^{\otimes s_1}) \otimes
                       H^0(K_X^{\otimes s_2})&
                       \longrightarrow & H^0(\theta^{\otimes s_1}) \otimes
                       H^0(\theta^{\otimes s_2}) \cr  @VV  V @VV
                       V\cr
                       H^0(K_X^{\otimes s_1+s_2})& \twoheadrightarrow &
                       H^0(\theta^{\otimes s_1+s_2})  \ ,
                       \endmatrix
                       $$
                       where the horizontal arrows are induced by restricting global
                       sections from
                       $X$
                       to $C$. The surjectivity of the
                       bottom horizontal arrow follows from  \smallemma \ns. 
                       Then it is clear that
                       the surjectivity of the left hand side vertical arrow implies the
                       surjectivity of
                       right hand side vertical arrow.

                        Now we prove 2). Assume  the maps  
$$H^0(\theta^{\otimes s'})
                        \otimes H^0(\theta^{\otimes s_2})
                        @> \alpha(s',s_2;1) >> H^0(\theta^{\otimes
                          s'+s_2})$$
 surject for
                        all
                        $1
                        \leq s' \leq s_1$.
                        We are going to show the surjectivity
                        of
                        $\alpha(s',s_2;m')$ for all $1 \leq m' \leq
                        m$ and  all $0 \leq s' \leq
                        s$. The proof is by induction on both on $m'$. If $m'=1$,
                        the result is our hypothesis, the case $s'=0$ being obvious. Now
                        we assume the result true for
                        $m'-1$ and we will prove it for $m'$ by induction on $s'$.
                        If $s'=0$, the result is
                        again obvious. Now we assume the result to be true for $s'-1$ and
                        we will prove it for $s'$. Consider the following  commutative
                        diagram:

                         $$\matrix
                         H^0(L_{m'}^{\otimes s'-1}) \otimes  H^0(L_{m'}^{\otimes s_2}) &
                         \hookrightarrow & H^0(L_{m'}^{\otimes s'}) \otimes
                         H^0(L_{m'}^{\otimes s_2})&
                         \twoheadrightarrow & H^0(L_{m'-1}^{\otimes s'}) \otimes
                         H^0(L_{m'}^{\otimes s_2}) \cr @VV   V @VV
                          V @VV  V\cr
                          H^0(L_{m'}^{\otimes s'+s_2-1}) & \hookrightarrow &
                          H^0(L_{m'}^{\otimes s'+s_2})& \twoheadrightarrow &
                          H^0(L_{m'-1}^{\otimes s'+s_2})  \ .
                          \endmatrix
                          $$

                           The last top and bottom horizonta
l arrows are surjective by
\smallemma \ns.  The left hand side vertical arrow is $\alpha(s'-1,s_2;
m')$ and surjects by induction hypothesis on $s'$. The right hand side
vertical arrow is the composition of the map
$$H^0(L_{m'-1}^{\otimes s'}) \otimes
H^0(L_{m'}^{\otimes s_2}) \longrightarrow H^0(L_{m'-1}^{\otimes
s'}) \otimes  H^0(L_{m'-1}^{\otimes s_2}) \ ,$$
 whose surjectivity
follows from \smallemma \ns, and the map of multiplication of global
sections on
$X_{m'-1}$, $\alpha{(s',s_2; m'-1)}$, which is surjective by
induction hypothesis on $m'$. Then it follows from chasing the
diagram that the middle vertical arrow
$\alpha(s',s_2; m')$ surjects.

Finally we prove 3)  by induction on  $m'$. If
 $m'=1$, there is nothing to prove. Now we assume the result
 true for $m'-1$ with $m' \geq 2 $ and we will prove it for $m'$.
 Consider the following  commutative diagram:

  $$\matrix
  H^0(L_{m'}^{\otimes s_1-1}) \otimes  H^0(L_{m'}^{\otimes s_2}) &
  \hookrightarrow & H^0(L_{m'}^{\otimes s_1}) \otimes
  H^0(L_{m'}^{\otimes s_2})&
  \twoheadrightarrow & H^0(L_{m'-1}^{\otimes s_1}) \otimes
  H^0(L_{m'}^{\otimes s_2}) \cr @VV   V @VV
   V @VV  V\cr
    H^0(L_{m'}^{\otimes s_1+s_2-1}) & \hookrightarrow &
    H^0(L_{m'}^{\otimes s_1+s_2})& \twoheadrightarrow &
    H^0(L_{m'-1}^{\otimes s_1+s_2})
    \cr @VV   V @VV  V
    @VV  V\cr
    0 &
    \longrightarrow & W_{m'} &
    = & W_{m'-1} \ . \cr
    \endmatrix
    $$

     The last top and bottom horizontal arrows are surjective by
     \smallemma \ns.
         The
     left hand side vertical arrow is $\alpha(s_1-1,s_2;m')$. This map
     surjects by part 2) of this lemma.  The right hand side
     first vertical arrow is the composition of the map
     $$H^0(L_{m'-1}^{\otimes s_1}) \otimes
     H^0(L_{m'}^{\otimes s_2}) \longrightarrow H^0(L_{m'-1}^{\otimes
     s_1})
     \otimes  H^0(L_{m'-1}^{\otimes s_2}) \ ,$$ 
and the map of
     multiplication of global sections on $X_{m'-1}$,
     $\alpha(s_1,s_2;{m'-1})$. The first of the two above mentioned maps
     is surjective by \smallemma \ns. Therefore the image of the right
     hand side first vertical arrow is equal to the image of
     $\alpha(s_1,s_2;{m'-1})$. By induction hypothesis, the image of
     $\alpha(s_1,s_2;{m'-1})$ is the inverse image of
      the image of $\alpha(s_1,s_2; 1)$ 
     by the natural map of restriction of global
     sections
     $$H^0(X_{m'-1},L_{m'-1}^{\otimes s_1+s_2}) \longrightarrow
     H^0(C,\theta^{\otimes s_1+s_2})\ .$$ By chasing the diagram it
     follows that the image of
     $\alpha(s_1,s_2,{m'})$  is the inverse image of
      the image of $\alpha(s_1,s_2,{m'-1})$ by the restriction
     map
     $$H^0(X_{m'},L_{m'}^{\otimes s_1+s_2})
     \longrightarrow H^0(X_{m'-1},L_{m'-1}^{\otimes s_1+s_2})\ .$$
     So finally the image of $\alpha(s_1,s_2;{m'})$ is the inverse
     image of
      the image of $\alpha(s_1,s_2; 1)$ by the restriction 
     map $$H^0(X_{m'},L_{m'}^{\otimes
     s_1+s_2})
     \longrightarrow H^0(C,\theta^{\otimes s_1+s_2})\ .$$
     Note finally
     that the isomorphism between $W_{m'}$ and
     $W_{m'-1}$ at each step of the inductive process tells us that the
     codimensions of  the images of $\alpha(s_1,s_2;{m'})$ in
     $H^0(X_{m'},L_{m'}^{\otimes s_1+s_2})$ are all equal for all $1 \leq
     m' \leq m$.
     \qed

     As we have previously mentioned in order to prove  \main  we need
      to study certain maps of multiplication of global sections of line
      bundles on the curve $C$ defined in  \notation \ns.  We will prove
      \main first when the dimension of $m$ is odd. On the curve $C$ we need
      the following 

      \proclaim{\oddcurve} Let $C$ be a smooth curve, let $\theta$ be
       an ample and
       base-point-free line bundle on $C$
       and let $m$ be an
       odd
       natural number such that $\theta^{\otimes m}=K_C$. Let $\varphi: X
       \longrightarrow Y$ be the morphism induced by a subseries $W'$ of
       $|\theta|$ without base points, and assume that
       the degree of $\varphi$ is $3$ and that $Y$ is a (smooth) rational
       normal curve of degree $r$. Then

         \itemitem{1.1)}$$\varphi_*(\Cal O_C)=\Cal O_{\bold P^1}
         \oplus \Cal O_{\bold P^1}(-\frac {mr+2}{2}) \oplus \Cal O_{\bold
         P^1}(-mr-2)) \ .$$ In particular,  $r$ is
         even.

          \itemitem{1.2)} If $m \geq 3$, then
          $\varphi$ is in fact induced by the complete linear series
            of $\theta$.

             \itemitem{2.1)} The map $\beta_n$
              in $H^0(\theta^{\otimes n+1})$ surjects if and only if
              $n \neq
              \frac{m-1}{2}, m$;

              \itemitem{2.2)} The map
               $\beta(\frac{m+1}{2},l)$ surjects
               in $H^0(\theta^{\otimes\frac{m+1}{2}+l})$ if $0 \leq l \leq
               \frac{m+1}{2}$.
              
                \itemitem{2.3)} The image of
                $\beta(s_1,s_2)$ is the same subspace
                $U'$ of
                $H^0(\theta^{\otimes \frac{m+1}{2}})$ for all

                $s_1,s_2$ such that $s_1+s_2=\frac{m+1}{2}, s_1, s_2 \geq 1$. The
                codimension of
                $U'$ is
                $\frac{r}{2}$.

                \itemitem{2.4)} The map $\beta(s_1,s_2)$ surjects if
                $s_1, s_2 \geq 0$ and $s_1+s_2 \leq \frac{m-1}{2}$.

                   \item{3)} The ring $R_\theta=\bigoplus_{s=1}^\infty (R_\theta)_s$,
                     where $(R_\theta)_s=H^0(\theta^{\otimes s})$, is
                     generated by its part of
                     degree

                     $1$ and
                     $\frac{r}{2}$ generators in degree $\frac{m+1}{2}$.

                     \endproclaim

                      \noindent{(2.5) \it Proof of \oddcurve \ns.}
                       We first prove 1.1).
                      Since $\varphi$ is finite and $Y$ is smooth, $\varphi$ is finite and
                       flat, then $\varphi_*(\Cal O_C)$ is a vector bundle of rank $3$.
                       The bundle $\varphi_*(\Cal O_C)$ splits as direct sum of line bundles
                       over $\bold P^1$. Let

                       $$\varphi_*(\Cal O_C)=\Cal O_{\bold P^1}
                       \oplus \Cal O_{\bold P^1}(-a_1) \oplus \Cal O_{\bold P^1}(-a_2)$$
                       for integers $a_1 \leq a_2$, which are strictly positive because $C$ is
                      connected, hence $h^0(\Cal O_C)=1$. Recall that
                       $K_C=\theta^{\otimes m}$ and that
                       $\theta=\varphi^*(\Cal O_{\bold P^1}(r))$. Then, on the one hand,
                       $$\varphi_*(K_C)=\Cal O_{\bold P^1}(mr) \oplus \Cal O_{\bold P^1}(mr-a_1)
                       \oplus   \Cal O_{\bold P^1}(mr-a_2) \ ,$$
by projection formula. On the other
                       hand, $$\varphi_*(K_C)=(\varphi_*(\Cal O_C))^* \otimes K_{\bold P^1} = \Cal
                       O_{\bold P^1}(a_2-2)
                       \oplus
                       \Cal O_{\bold P^1}(a_1-2)
                       \oplus   \Cal O_{\bold P^1}(-2) \ ,$$ by relative duality. Then
                       $$\displaylines{ mr=a_2-2 \cr
                       mr-a_1=a_1-2 \cr
                       mr-a_2=-2 \ ,}$$
                       from where $a_1=\frac {mr+2}{2}$ and $a_2=mr+2$. Moreover,
                       since $m$ is odd and $mr+2$ has to be even, $r$ is
                       even. This completes the proof of 1.1).

                       \medskip

                       Now we prove 1.2). It follows from 1.1) and projection formula that

                          $$\displaylines{ H^0(\theta) = H^0(\varphi_*(\theta))=H^0(\Cal
                          O_{\bold P^1}(r)) \cr
                          \oplus \Cal O_{\bold P^1}(\frac {(2-m)r-2}{2}) \oplus \Cal O_{\bold
                          P^1}((1-m)r-2))) \  .} $$

                          Then $H^0(\varphi_*(\theta))=H^0(\Cal O_{\bold P^1}(r))$,
                           for $m \geq 3$. Therefore $\varphi$ is
                           induced by the complete $|\theta|$.

                           We have seen that $\varphi_*(\Cal O_C)$ and $\Cal O_{\bold P^1}
                            \oplus \Cal O_{\bold P^1}(-\frac {mr+2}{2}) \oplus \Cal O_{\bold
                            P^1}(-mr-2))$
                              are isomorphic as sheaves of
                              modules over $\Cal O_{\bold P^1}$. We call $E=\Cal O_{\bold
                                P^1}(-\frac {mr+2}{2}) \oplus \Cal O_{\bold
                                P^1}(-mr-2))$. Since $\varphi_*(\Cal O_C)$ is also a
                                sheaf of algebras over $\Cal O_{\bold P^1}$, then the algebra
                                structure of $\varphi_*(\Cal O_C)$ makes $\Cal O_{\bold P^1}
                                \oplus \Cal O_{\bold P^1}(-\frac {mr+2}{2}) \oplus \Cal O_{\bold
                                P^1}(-mr-2))$  into
                                a sheaf of $\Cal O_{\bold P^1}$-algebras whose multiplication
                                decomposes like
                                this (see \Miranda \ns, \HM or \GP \ns):
                                The map $$\Cal O_{\bold P^1} \otimes \Cal O_{\bold P^1}
                                  \longrightarrow \Cal O_{\bold P^1}$$
                                  is the ring multiplication in $\Cal O_{\bold P^1}$.
                                  The maps
                                  $$\matrix \Cal O_{\bold P^1} \otimes  \Cal O_{\bold P^1}(-\frac
                                    {mr+2}{2}) &
                                    \longrightarrow &  \Cal O_{\bold P^1}(-\frac {mr+2}{2}) \cr
                                    \Cal O_{\bold P^1} \otimes  \Cal O_{\bold P^1}(-mr-2) &
                                    \longrightarrow &  \Cal O_{\bold P^1}(-mr+2) \cr
                                      \Cal O_{\bold P^1}(-\frac
                                        {mr+2}{2}) \otimes \Cal O_{\bold P^1} &
                                        \longrightarrow  & \Cal O_{\bold P^1}(-\frac {mr+2}{2}) & \text{ and }\cr
                                         \Cal O_{\bold P^1}(-mr-2) \otimes \Cal O_{\bold P^1} &
                                         \longrightarrow  & \Cal
                                         O_{\bold P^1}(-mr+2) \endmatrix$$
                                         are the left and right module multiplication of $\Cal O_{\bold
                                           P^1}(-\frac {mr+2}{2})$ and $\Cal O_{\bold P^1}(-mr+2)$.
                                           Finally there is a map
                                           $$ E \otimes E \longrightarrow \varphi_*\Cal
 O_C $$
 whose structure depends on the cover $\varphi$ itself.
 We will use all this to
 prove 2).
Since $\theta=\varphi^*\Cal O_{\bold P^1}(r)$, it follows from 1) and
projection formula that 
$$\displaylines{H^0(\theta^{\otimes n})=H^0(\Cal O_{\bold P^1}(nr))
\oplus H^0(\Cal O_{\bold P^1}(\frac{(2n-m)r-2)}{2})) \oplus
H^0(\Cal O_{\bold P^1}((n-m)r-2))\ .}$$ 
Let us denote 
$$\matrix A(n)&=&H^0(\Cal O_{\bold P^1}(nr)) \cr
B(n)&=&H^0(\Cal O_{\bold P^1}(\frac{(2n-m)r-2)}{2})) \cr
C(n)&=&H^0(\Cal O_{\bold P^1}((n-m)r-2)) &\ .\endmatrix $$
The way in which global sections of $H^0(\theta^{\otimes s_1})$ and
global sections of $H^0(\theta^{\otimes s_2})$ multiply is governed by the ring
multiplication in $$\varphi_*\Cal O_C=\Cal O_{\bold P^1}
\oplus \Cal O_{\bold P^1}(-\frac {mr+2}{2}) \oplus \Cal O_{\bold
P^1}(-mr-2)) \ .$$ Thus
the map $\beta(s_1,s_2)$  splits as direct sum of the following maps
$$\displaylines{A(s_1) \otimes A(s_2) \longrightarrow
A(s_1+s_2) \cr
A(s_1) \otimes B(s_2)
\longrightarrow  B(s_1+s_2) \cr
A(s_1) \otimes C(s_2)
\longrightarrow C(s_1+s_2) \cr
B(s_1) \otimes A(s_2) \longrightarrow B(s_1+s_2) \cr
C(s_1) \otimes A(s_2) \longrightarrow C(s_1+s_2) \cr
[B(s_1) \oplus C(s_1)] \otimes [B(s_2) \oplus C(s_2)]  
\longrightarrow \cr
A(s_1+s_2) \oplus B(s_1+s_2) \oplus C(s_1+s_2) \oddformulaone .}$$

The first map $A(s_1) \otimes A(s_2) \longrightarrow
 A(s_1+s_2)$ is induced by the ring multiplication in $\Cal O_{\bold
   P^1}$ and it is therefore a map of multiplication of global sections
   of line bundles on $\bold P^1$. The maps second to fifth are induced
   by the left and right module structure of $E$ over $\Cal O_{\bold
     P^1}$ and are also multiplcation maps of global sections of line
     bundles on $\bold P^1$. The structure of the last map is more complex
     and depends on the structure of $\varphi$. Thus the image of
     $\beta(s_1,s_2)$ will be the sum of the images of those six maps. We
assume 
$s_1, s_2 \geq 0$. Then the first map is always nonzero and surjective,
hence $A(s_1+s_2)$ is always contained in the image of
$\beta(s_1,s_2)$. On the other hand the vanishing of the groups
$B(s_1), B(s_2), B(s_1+s_2)$ and $C(s_1), C(s_2), C(s_1+s_2)$ depends 
on the values of $s_1$ and $s_2$. If $B(s_1)$ does not vanish, that
is, if $s_1 \geq \frac{m+1}{2}$, then $B(s_1+s_2)$ does not vanish
either and the map
$$B(s_1) \otimes A(s_2) \longrightarrow B(s_1+s_2)$$ is
surjective. Therefore in this case $B(s_1+s_2)$ 
is contained also in the image of $\beta(s_1,s_2)$. 
If $C(s_1)$ does not vanish, that is, if $s_1 \geq m+1$, then
$C(s_1+s_2)$ does not vanish 
either and the map $$C(s_1) \otimes A(s_2) \longrightarrow
C(s_1+s_2)$$ is surjective. Therefore in this case $C(s_1+s_2)$ is
contained  in the image of $\beta(s_1,s_2)$.
There are analogous statements for $B(s_2)$ and $C(s_2)$.

  \medskip

   With this information we prove now 2.1), 2.2), 2.3) and 2.4).
   We prove first 2.1). The multiplication map $\beta_n=\beta(n,1)$
    splits as shown in \oddformulaone \ns.
    If $0 \leq n < \frac{m-1}{2}$, $B(n), C(n), B(n+1)$ and $C(n+1)$ all
    vanish, therefore $H^0(\theta^{\otimes n+1})=A(n+1)$, and the map
    $$H^0(\theta^{\otimes n}) \otimes H^0(\theta) @> \beta_n >>
    H^0(\theta^{\otimes n+1}) $$ is
    the map
    $$A(n) \otimes A(1) \longrightarrow A(n+1) \ ,$$
    which is surjective.
    On the other hand,
    $B(\frac{m-1}{2})=C(\frac{m-1}{2})=C(\frac{m+1}{2})=0$ but
    $B(\frac{m+1}{2})$ has dimension $\frac r 2 > 0$. Then
$H^0(\theta)=A(1)$,
$H^0(\theta^{\otimes\frac{m-1}{2}})=A(\frac{m-1}{2})$ and
$H^0(\theta^{\otimes\frac{m+1}{2}})=A(\frac{m+1}{2}) \oplus
B(\frac{m+1}{2})$. Then
$$H^0(\theta^{\otimes \frac{m-1}{2}}) \otimes H^0(\theta) @>
\beta_{\frac{m-1}{2}} >> H^0(\theta^{\otimes\frac{m+1}{2}})$$ is
the map $$A(\frac{m-1}{2}) \otimes A(1) \longrightarrow
A(\frac{m+1}{2}) \ , $$
whose image has codimension $\frac r 2$ in
$H^0(\theta^{\otimes\frac{m+1}{2}})$.
If $\frac{m+1}{2} \leq n \leq m-1$, then
$H^0(\theta^{\otimes n}) = A(n) \oplus B(n)$ and
$H^0(\theta^{\otimes n+1}) = A(n+1) \oplus B(n+1)$, with
$A(n), B(n)$, $A(n+1)$ and $B(n+1)$ all different from $0$. Then
according to (\oddformulaone \ns) if $\frac{m+1}{2} \leq n \leq m-1$ the map
$\beta_n$ decomposes as direct sum of
$$\displaylines{A(n) \otimes A(1) \longrightarrow A(n+1) \text{ and } \cr
B(n) \otimes A(1) \longrightarrow B(n+1) \ ,}$$
which are surjective maps, therefore the image of $\beta_n$ is $A(n+1)
\oplus B(n+1)=H^0(\theta^{\otimes n+1})$.
On the other hand $C(m)=0$ but $C(m+1)$ has dimension $r-1 > 0$. Thus
$H^0(\theta^{\otimes m})=A(m) \oplus B(m)$ but
$$H^0(\theta^{\otimes m+1})=A(m+1) \oplus B(m+1) \oplus C(m+1) \ .$$
The map $\beta_m$ decomposes as direct sum of
$$\displaylines{A(m) \otimes A(1) \longrightarrow A(m+1) \text{ and } \cr
B(m) \otimes A(1) \longrightarrow B(m+1) \ ,}$$
and then its image is $A(m+1) \oplus B(m+1)$, which has codimension
$r-1$ in $H^0(\theta^{\otimes m+1})$. 
Finally if $n \geq m+1$, $H^0(\theta^n) = A(n) \oplus B(n) \oplus
C(n)$, and $A(n), B(n)$ and $C(n)$ are all nonzero. Then $\beta_n$
decomposes as direct sum of
$$\displaylines{A(n) \otimes A(1) \longrightarrow A(n+1) \cr
  B(n) \otimes A(1) \longrightarrow B(n+1) \cr
C(n) \otimes A(1) \longrightarrow C(n+1) \ ,}$$
and its image is $A(n+1) \oplus B(n+1) \oplus
C(n+1)=H^0(\theta^{n+1})$. This completes the proof of 2.1).

Now we prove 2.2).
Recall that if $\frac{m+1}{2} \leq n \leq m$, then 
$H^0(\theta^{\otimes n}) =A(n) \oplus B(n)$ and that, 
if $0 \leq l \leq \frac{m-1}{2}$, then
$H^0(\theta^{\otimes l})=A(l)$. 
As we argued in the proof of 2.1, the map $\beta(\frac{m+1}{2},l)$ is
governed by the algebra structure of $\varphi_*\Cal O_C$. 
Then, if $0 \leq l \leq \frac{m-1}{2}$,  $\beta(\frac{m+1}{2},l)$
decomposes as direct sum of
$$\displaylines{A(\frac{m+1}{2}) \otimes A(l) \longrightarrow
  A(\frac{m+1}{2}+l) \cr
B(\frac{m+1}{2}) \otimes A(l) \longrightarrow B(\frac{m+1}{2} +l) \
.}$$
 Then the image of $\beta(\frac{m+1}{2},l)$ is $ A(\frac{m+1}{2}+l)
\oplus B(\frac{m+1}{2} +l) = H^0(\theta^{\otimes \frac{m+1}{2}+l})$. 

On the other
hand, $H^0(\theta^{\otimes \frac{m+1}{2}})= A(\frac{m+1}{2}) \oplus
B(\frac{m+1}{2}) $.
Therefore the map
$\beta(\frac{m+1}{2}, \frac{m+1}{2})$ decomposes as direct sum of
$$\displaylines{A(\frac{m+1}{2}) \otimes A(\frac{m+1}{2}) \longrightarrow
  A(m+1) \cr
  B(\frac{m+1}{2}) \otimes A(\frac{m+1}{2}) \longrightarrow B(m+1) \cr
  A(\frac{m+1}{2}) \otimes B(\frac{m+1}{2}) \longrightarrow B(m+1) \cr
  B(\frac{m+1}{2}) \otimes B(\frac{m+1}{2}) @> \eta >> \cr 
A(m+1)
  \oplus B(m+1) \oplus C(m+1)  \ \ \oddformulatwo .}$$
  Since the first three maps are surjective
  We know that $A(m+1) \oplus B(m+1)$ is
  contained in the image of $\beta(\frac{m+1}{2},\frac{m+1}{2})$.
  Now we need to look more carefully at the last one of the above
  maps.
  For this  we have to study more closely the
  ring structure of  $\varphi_*(\Cal O_C)$.  Let
  $$ \varphi_*(\Cal O_C) \otimes \varphi_*(\Cal O_C) @> \mu >>
  \varphi_*(\Cal O_C) $$ be the map defined by the ring
  multiplication in
 $\varphi_*(\Cal O_C)$.  Then the projection to $\Cal O_{\bold
P^1}(-mr-2)$ of the restriction of $\mu$ to
$\Cal O_{\bold P^1}(\frac{-mr-2}{2}) \otimes \Cal
O_{\bold P^1}(\frac{-mr-2}{2}))$ is nonzero. Otherwise there
would be an integral
$\Cal O_{\bold P^1}$-subalgebra of $\varphi_*(\Cal O_C)$ of rank $2$,
namely $\Cal O_{\bold P^1} \oplus \Cal O_{\bold
  P^1}(\frac{-mr-2}{2})$. This is impossible because the rank of an
  integral subalgebra of $\varphi_*\Cal O_C$ must divide the rank of
  $\varphi_*\Cal O_C$.
  Now if the projection to $\Cal O_{\bold
  P^1}(-mr-2)$ of the restriction of $\mu$ to
  $\Cal O_{\bold P^1}(\frac{-mr-2}{2}) \otimes \Cal
  O_{\bold P^1}(\frac{-mr-2}{2}))$ is nonzero, it is an
  isomorphism. Then the projection of $\eta$ to $C(m+1)$ is
  surjective. Since the image of the three first maps in
  (\oddformulatwo \ns) is $A(m+1)
  \oplus B(m+1)$ we may conclude that the image of $\beta(\frac{m+1}{2},
  \frac{m+1}{2})$ is
  $$A(m+1) \oplus B(m+1) \oplus C(m+1) = H^0(\theta^{\otimes m+1}) \ .$$

We proceed now to prove 2.3). 
Let $s_1, s_2$ satisfy 
$s_1+s_2=\frac{m+1}{2}$ and $s_1,s_2 \geq 1$. Recall that
$$H^0(\theta^{\otimes \frac{m+1}{2}})=
A(\frac{m+1}{2}) 
\oplus B(\frac{m+1}{2}) \ .$$
However, if
$s_1+s_2=\frac{m+1}{2}$ and $s_1,s_2\geq 1$, then $s_1,s_2 \leq
\frac{m-1}{2}$, hence
$H^0(\theta^{\otimes s_1})=A(s_1)$ and
$H^0(\theta^{\otimes s_2})=A(s_2)$. Thus the map $\beta(s_1,s_2)$ is 
$$A(s_1) \otimes A(s_2) \longrightarrow A(\frac{m+1}{2})\ . $$ 
Thus the image
of $\beta(s_1,s_2)$ is, for all
$s_1+s_2=\frac{m+1}{2}$, $s_1,s_2 \geq 1$, equal to
$A(\frac{m+1}{2})=H^0(\Cal O_{\bold P^1}(\frac{m+1}{2}r))$ 
which has codimension
$\frac {r}{2}$ in $H^0(\theta^{\otimes \frac{m+1}{2}})$. 

Finally we prove 2.4). If $0 \leq s \leq \frac{m-1}{2}$
then $H^0(\theta^{\otimes s})=A(s)$. Then if $s_1,s_2 \geq 0$ and
$s_1+s_2 \leq \frac{m-1}{2}$, the map $\beta(s_1,s_2)$ is
the multiplication map 
$$A(s_1) \otimes A(s_2) \longrightarrow A(s_1+s_2)$$ which is
surjective.

Now we prove 3). Note that the vector subspace of $(R_\theta)_n$
generated by the elements of degree smaller than $n$ is the sum of
the images of the maps $\beta(s_1,s_2)$ with $s_1, s_2 \geq 1$ and
$s_1+s_2=n$. We saw in 2.1) that the maps $\beta_l$ surject if $l <
\frac{m-1}{2}$. Then the parts of $R_\theta$ of degree less than
$\frac{m+1}{2}$ are generated by the elements of degree $1$. 2.3)
tells us that the elements of degree less than $\frac{m+1}{2}$
generates a subspace of $(R_\theta)_{\frac{m+1}{2}}$ 
of  codimension $\frac r 2$, 
henceforth  we need $\frac r 2$ linearly independent 
elements of degree $\frac{m+1}{2}$ 
to generate $(R_\theta)_{\frac {m+1}{2}}$.
If $\frac{m+1}{2} \leq l \leq m-1$, we saw also in 2.1)
that the maps $\beta_l$ surject, therefore the parts of degree
$\frac{m+1}{2}, \dots, m$  are generated by the elements of $R_\theta$ of
degree less than or equal to $\frac{m+1}{2}$. From 2.2) it follows
that the map $\beta(\frac{m+1}{2},\frac{m+1}{2})$ surjects, therefore
$(R_\theta)_{m+1}$ is also generated by the elements of
$(R_\theta)_{\frac{m+1}{2}}$. Finally, if $l \geq m+1$ we saw in 2.1) that the
map $\beta_l$ surjects, thus the parts of $R_\theta$ of degree greater than
or equal to $m+2$ are generated by the elements of degree
$m+1$. Summarizing the ring $R_\theta$ is generated by its elements of degree
$1$ and by $\frac {r}{2}$ elements of degree $\frac{m+1}{2}$.
\qed

 \bigskip

 Now we can prove \main when the dimension of the variety of general
  type $X$ is odd. Precisely we prove

  \proclaim{\oddmain}  Let $X$ be a smooth,  pluriregular 
   variety of general type
   of odd dimension $m \geq 3$  
    and let $Y
   \subset
   \bold P^{r+m-1}$ be an $m$-fold of minimal degree  $r$ . Assume that
   $K_X$ is base-point-free and let
   $\pi: X \longrightarrow Y$ be the canonical morphism of $X$.
      Assume that $\pi$ is generically finite morphism of degree $3$. Then
   the canonical ring of $X$ is generated by its part of
   degree
$1$ and
$\frac{r}{2}$ generators in degree $\frac{m+1}{2}$. 
\endproclaim

\noindent {\it Proof.}
Let $C$ be the intersection of $m-1$ general members  of
$|K_X|$. Then $C$ is a smooth and irreducible curve in $X$.
Let $\theta=K_X|_C$. By adjunction
$\theta^{\otimes m} = K_C$. We set $\varphi=\pi|_C$. Since
$C$ has been taken  general, we may, and actually do, choose
$C$ so that $\varphi$ is finite onto a smooth rational normal curve $D$
of degree $r$ in $\bold P^r$.

\medskip

 As we have seen in the proof of \oddcurve \ns, in order to find out
  the degrees of the generators of $R$, we only need to study the images of
  the maps $\alpha(s,t)$ for all $s,t \geq 1$.  In fact, since many
  of them will happen to be surjective it will suffice, for our purpose,
  to study only some of them.
  More precisely the statement of the theorem is true if
  the following happens:

  \medskip
   \item{1)} 
   $H^0(K_X^{\otimes n}) \otimes H^0(K_X) @> \alpha_n
>> H^0(K_X^{\otimes n+1})$  surject for all $n \geq 0$,  
$n \neq \frac{m-1}{2}, m$;  
\item{2)}
$H^0(K_X^{\otimes \frac{m+1}{2}}) \otimes H^0(K_X^{\otimes
\frac{m+1}{2}})
@> \alpha(\frac{m+1}{2},\frac{m+1}{2}) >> H^0(K_X^{\otimes m+1})$
surjects; 
\item{3)} the images  of \newline 
\centerline{$H^0(K_X^{\otimes
s_1}) \otimes H^0(K_X^{\otimes s_2}) @>\alpha(s_1,s_2) >>
H^0(K_X^{\otimes
\frac{m+1}{2}})$} \newline
for all $s_1, s_2 \geq 1, s_1+s_2=\frac{m+1}{2}$
are equal of codimension $\frac{r}{2}$.

Recall that $\varphi=\pi|_C$ is induced by the complete series of
 $\theta$ and that $C$, $\theta$ and $\varphi$ satisfy the hypothesis
 of \oddcurve \ns. Now to  prove 1), 2) and 3) we will relate the maps
 appearing in 1), 2) and 3)
 with the multiplication maps on
 $C$ studied in \oddcurve \ns.
 To establish the relation between them we will use among other things
 \superlemma \ns.  
1) and 2) follows directly from \superlemma (2)  and \oddcurve
  (2.1) and 2.4 (2.2).   We now prove 3). By \oddcurve (2.4) we know
  that  the maps $\beta(s_1,s_2)$ are surjective if $s_1 + s_2 \leq
  \frac{m-1}{2}$. Let us know fix $s_1, s_2 \geq 1$ such that $s_1 + s_2
  = \frac{m+1}{2}$. Then by \superlemma (3), we know that  the codimension of
the image of 
$\alpha(s_1,s_2)$ in $H^0(K_X^{\otimes \frac{m+1}{2}})$ is the same as
the codimension of the image of
$\beta(s_1,s_2)$ in
$H^0(\theta^{\otimes \frac{m+1}{2}})$. We also know by \superlemma (3) that
the image of $\alpha(s_1,s_2)$ is  the inverse image
by the obvious map of restriction of
sections
$$ H^0(K_X^{\otimes \frac{m+1}{2}}) \longrightarrow
H^0(\theta^{\otimes \frac{m+1}{2}}) $$
of
the image of $\beta(s_1,s_2)$. 
 Finally \oddcurve (2.3) tells
 the image of $\beta(s_1,s_2)$ in
 $H^0(\theta^{\otimes \frac{m+1}{2}})$ are the
 same for all $s_1, s_2 \geq 1$ such that $s_1+s_2=\frac{m+1}{2}$ and
 have codimension $\frac r 2$. Then all this implies that the images of
 $\alpha(s_1,s_2)$ are equal and of codimension $\frac r 2$ in
 $H^0(K_X^{\otimes \frac{m+1}{2}})$ for all $s_1, s_2 \geq 1$ such that
 $s_1+s_2=\frac{m+1}{2}$.
  \qed

        \proclaim{\oddremark} \oddcurve states
        explicitely that the ring  
$\bigoplus_{n \geq
        0}H^0(C,\theta^{\otimes n})$ has in each degree the same number of
        generators as the canonical ring of $X$, $\bigoplus_{n \geq
        0}H^0(X,K_X^{\otimes n})$. The way \oddmain is proved shows that
        the same is true also for the ``intermediate" rings
        $\bigoplus_{n
        \geq 0}H^0(X_{m'},L_{m'}^{\otimes n})$, $2 \leq m' \leq m-1$.
        \endproclaim

         Knowing the generators of the canonical ring of $X$ tells us about the
         generators of the pluricanonical rings of $X$ and in particular, about
         which pulricanonical rings are generated in degree $1$. We state te
         result we obtain in this regard in the next Theorem. Before that we
         recall the following definition due to M. Green:

         \proclaim{\defNo} Let $L$ be a line
          bundle on a projective variety $X$. We say that $L$ satisfies property
          $N_0$ if
          $L$ is very ample and the image of $X$ by the embedding induced by
          $|L|$ is projectively normal. 
          \endproclaim

          \proclaim{\oddNo} Let $X$ be a smooth, pluriregular
           variety of general type
           of odd dimension $m \geq 3$
           and let $Y
           \subset
           \bold P^{r+m-1}$ be an $m$-fold of minimal degree $r$. Assume
           that
           $K_X$ is ample base-point-free and let
           $\pi: X \longrightarrow Y$ be the canonical morphism of $X$.
           Assume that $\pi$ has degree $3$.  Then the line
           bundle
           $K_X^{\otimes n}$ satisfies property $N_0$
           if and only if $n \geq \frac{m+1}{2}$.
           \endproclaim

            \noindent {\it Proof.}  Let
             $C$ be as in \oddmain \ns. We carry out the proof in several steps.

             {\it Step 1.} $K_X^{\otimes n}$ does not satisfy property $N_0$ if
              $n <
              \frac{m+1}{2}$.

               Indeed, if
               $n <
               \frac{m+1}{2}$, then
               $$H^0(K_X^{\otimes n}) \longrightarrow H^0(K_X^{\otimes
                 n}|_C) $$ surjects by
               \smallemma and
               $$ H^0(K_X^{\otimes n}|_C)=H^0(\theta^{\otimes n})=H^0(\Cal O_{\bold
               P^1}(nr)) \ ,$$
                hence the $n$-pluricanonical morphism
                of $X$ maps $C$ to a rational normal curve, therefore $K_X^{\otimes
                  n}$ is not very ample and does
                  not satisfy property $N_0$. 

                  \medskip

                   Now to see that $K_X^{\otimes n}$
                    satisfies property $N_0$ for given $n \geq
                    \frac{m+1}{2}$,   
  it suffices  to prove the surjectivity of the
                    following maps:

                    $$H^0(K_X^{\otimes ln}) \otimes H^0(K_X^{\otimes n}) @>
                     \alpha(ln,n) >> H^0(K_X^{\otimes (l+1)n}), \text{ for all } l \geq
                     1 \ .$$

                     \medskip

                      {\it Step 2.} $K_X^{\otimes n}$ satisfies property $N_0$ if
                       $n =
                       \frac{m+1}{2}$.

                       The map $\alpha(\frac{m+1}{2}, \frac{m+1}{2})$ surjects as shown 
                        in
                         2) of the proof of \text{\oddmain \ns.}   On the other hand, the
                         surjectivity of
                         $\alpha(ln,n)$, $l \geq 2$ and $n = \frac{m+1}{2}$ follows from the
                         surjectivity of $\alpha_{n'}$ for all $n' \geq m+1$, which was also
                         shown in  1)  of the proof of \oddmain \ns.

                          \medskip

                          {\it Step 3.}  $K_X^{\otimes n}$ satisfies property $N_0$ if
                           $
                           \frac{m+1}{2} < n \leq m$.
We first show that  
$\alpha(n,n)$
surjects. Since $\alpha(\frac{m+1}{2}, \frac{m+1}{2})$ surjects and so
do the maps $\alpha_{n'}$ if $n' \geq m+1$, the following multiplication map
$$H^0(K^{\otimes \frac{m+1}{2}}) \otimes H^0(K_X^{\otimes
  (n-\frac{m+1}{2})}) \otimes H^0(K^{\otimes \frac{m+1}{2}}) \otimes
H^0(K_X^{\otimes 
  (n-\frac{m+1}{2})}) @> \gamma >> H^0(K_X^{\otimes 2n}) $$
also surjects. On the other hand $\gamma$ factorizes through
$\alpha(n,n)$, which is therefore surjective. 

Finally the surjectivity of the maps $\alpha(ln,n)$ for all $l \geq 2$ and $
\frac{m+1}{2} < n \leq m$ follows from the surjectivity of
$\alpha_{n'}$ for all $n' \geq m+1$. 

{\it Step 4.} $K_X^{\otimes n}$ satisfies property $N_0$ if
$n \geq m+1$. The surjectivity of the maps $\alpha_{n'}$ for all $n'
\geq m+1$ implies the surjectivity of the maps $\alpha(ln,n)$ for all
$l \geq 1$ and $n \geq m+1$.
 \qed

\bigskip

 Now we proceed to prove \main when the dimension of $X$ is even. 
As we did in the
odd dimensional case, we need study first multiplication maps on a
curve $C$, which is a complete intersection on $X$.  
Thus we prove a result analogous to \oddcurve \ns.

\proclaim{\evencurve} Let $C$ be a smooth curve, let $\theta$ be
an ample and 
base-point-free line bundle on $C$ 
and let $m$ be an
even 
natural number such that $\theta^{\otimes m}=K_C$. Let $\varphi: C
\longrightarrow Y$ be the morphism induced by $|\theta|$ and assume that
the degree of $\varphi$ is $3$ and that $Y$ is a (smooth) rational
normal curve of degree $r$. Then  
 
\itemitem{1.1)}$$\varphi_*(\Cal O_C)=\Cal O_{\bold P^1}
\oplus \Cal O_{\bold P^1}(-\frac {mr+2}{2}) \oplus \Cal O_{\bold
P^1}(-mr-2)) \ .$$

 \itemitem{1.2)} The map
 $\varphi$ is induced by the complete linear series
   of $\theta$.

   \medskip

     \itemitem{2.1)} The map $\beta_n$ surjects if $n \neq \frac{m}{2},$ 
      $m,
        m+1$; if $n=m+1$ and $r \geq 2$; and if $n=m$ and $r=1$.

        \itemitem{2.2)} The images of the maps $\beta(s_1,s_2)$ are equal and of
           codimension $r$  in \linebreak 
$H^0(\theta^{\otimes \frac {m+2}{2}})$ if $s_1, s_2
               \geq 1$ and $s_1 + s_2 = \frac{m+2}{2}$.

               \itemitem{2.3)} The images of the maps $\beta(s_1,s_2)$ are equal and of
                  codimension $r-1$  in $H^0(\theta^{\otimes m+1})$ if $s_1, s_2
                    \geq 1$ and $s_1 + s_2 = m+1$.

\itemitem{2.4)} The map $\beta(\frac{m+2}{2},\frac{m+2}{2})$ surjects. 

\itemitem{2.5)} The map $\beta(s_1,s_2)$ surjects if $s_1, s_2 \geq 0$
and $s_1+s_2 \leq \frac{m}{2}$.

\itemitem{2.6)} The map $\beta(s_1,s_2)$ surjects if $s_1, s_2 \geq
0$, $s_1 \geq \frac {m+2}{2}$, 
and $\frac {m+2}{2} \leq s_1+s_2 \leq m$. 

\medskip

\item{3)} The ring $R_\theta=\bigoplus_{s=1}^\infty (R_\theta)_s$,
  where $(R_\theta)_s=H^0(\theta^{\otimes s})$, is
  generated by its part of
  degree
  $r$
  generators in degree $\frac{m+2}{2}$ and $r-1$ generators in
  degree $m+1$.

  \endproclaim

   {\it Proof.} The proof of 1.1) and 1.2) is exactly like the proof of
    parts 1.1) and 1.2) of
    \oddcurve \ns.
    The proof of 2) is also like the proof of part 2) of \oddcurve \ns. Here
    we will only highlight the differences. As we saw in the proof of
    \oddcurve (2), the key to know the images of the maps
    $\beta(s_1,s_2)$ is the decomposition of $H^0(\theta^{\otimes s_1})$
    and $H^0(\theta^{\otimes s_2})$ as sum of blocks $A(s_1) \oplus B(s_1)
    \oplus C(s_1)$ and $A(s_2) \oplus B(s_2) \oplus C(s_2)$, i.e., the
    crucial information is to know which ones among $A(s_1), B(s_1),
    C(s_1), A(s_2), B(s_2), C(s_2)$ are zero. This information is obtained
    from 1).
    If $0 \leq n \leq \frac{m}{2}$, then $H^0(\theta^{\otimes
      n})=A(n)$. Thus, the map $\beta_n$ surjects if $0 \leq n \leq
\frac{m-2}{2}$. If $\frac{m+2}{2} \leq n \leq m$,  $H^0(\theta^{\otimes
  n})=A(n) \oplus B(n)$, hence the map $\beta_n$ also surjects if
  $\frac{m+2}{2} \leq n \leq m-1$. If $n \geq m+2$, $H^0(\theta^{\otimes
  n})=A(n) \oplus B(n) \oplus C(n)$, then the map $\beta_n$ surjects if
  $n \geq m+2$. 
If $r=1$, $H^0(\theta^{\otimes m+1}) = A(m+1) \oplus B(m+1)$, thus
$\beta_m$ surjects if $r=1$. 
If $r > 1$, then $H^0(\theta^{\otimes m+1}) = A(m+1) \oplus B(m+1)
\oplus C(m+1)$, therefore $\beta_{m+1}$ surjects if $r > 1$. 
This proves 2.1). Now we prove 2.2). Recall that if $0 \leq n \leq
\frac{m}{2}$, then $H^0(\theta^{\otimes 
  n})=A(n)$. On the other hand $H^0(\theta^{\otimes
  \frac{m+2}{2}})=A(\frac{m+2}{2}) \oplus B(\frac{m+2}{2})$. Since
$B(\frac{m+2}{2})= H^0(\Cal O_{\bold P^1}(r-1))$, then, if $s_1, s_2
\geq 0$ and $s_1 + s_2=\frac{m+2}{2}$ the image of the maps
$\beta(s_1,s_2)$ are equal to $A(\frac{m+2}{2})$, which has
codimension $r$ in $H^0(\theta^{\otimes \frac{m+2}{2}})$. 
To prove 2.3) recall that if $\frac{m+2}{2} \leq n \leq m$,  
$H^0(\theta^{\otimes
  n})=A(n) \oplus B(n)$ and that $H^0(\theta^{\otimes m+1}) = A(m+1)
    \oplus B(m+1)
    \oplus C(m+1)$.  Thus, if $s_1, s_2 \geq 0$ and $s_1+s_2=m+1$, then
    either $s_1 \geq \frac{m+1}{2}$ or $s_2 \geq \frac{m+1}{2}$. Let us
    say that $s_1 \geq \frac{m+1}{2}$, then $s_2 \leq \frac{m+1}{2}$.  Then the
    image of the maps $\beta(s_1,s_2)$ is $A(m+1) \oplus B(m+1)$. Since
    $C(m+1)=H^0(\Cal O_{\bold P^1}(r-2))$, the codimension of the image of
    $\beta(s_1,s_2)$ in $H^0(\theta^{\otimes m+1})$ is $r-1$.
    Part 2.4) follows from the fact that $H^0(\theta^{\otimes \frac
      {m+2}{2}})=A(\frac {m+2}{2}) \oplus B(\frac {m+2}{2})$ arguing as in
      the proof of \oddcurve (2.2).
      Part 2.5) follows from the fact that$H^0(\theta^{\otimes
        n})=A(n)$ if $0 \leq n \leq \frac m 2$.
        Now we prove 2.6). Recall that if $\frac{m+2}{2} \leq n \leq m$,
        then $H^0(\theta^{\otimes
          n})=A(n) \oplus B(n)$ and that if $0 \leq n \leq \frac{m}{2}$,
          then $H^0(\theta^{\otimes
            n})=A(n)$. This implies 2.6).

            Finally 3) follows from 2.1), 2.2), 2.3) and 2.4). \qed

             \bigskip
               Now
               from the previous \evencurve and using \smallemma and \superlemma
               we can prove \main for the even dimensional case: 

               \proclaim{\evenmain}  Let $X$ be a smooth, pluriregular variety
                of general type of even dimension $m$ with $K_X$ base-point-free.
                 Assume that the canonical morphism $\pi$ is
                generically finite of degree $3$ onto a variety $Y$ of minimal
                degree
                $r$. Then  the canonical ring of $X$ is generated by 
                its elements of degree $1$, $r$
                generators in  degree
                $\frac{m+2}{2}$ and $r-1$ generators in degree $m+1$.
                     \endproclaim

                \noindent{\it Sketch of proof.} The proof uses the same ideas as
                 the proof of \oddmain so we will just outline the key steps. As in
                 \oddmain we need to know the images of the maps $\alpha(s_1,s_2)$.
                 Precisely the
                 result follows from the following claims:

                 \medskip

                  \item{1)} The map $\alpha_n$ surjects if $1 \leq n \leq
                   \frac{m-2}{2}$.

                   \item{2)} The images of the maps $\alpha(s,t)$, when
$s+t=\frac{m+2}{2}$ and $s,t\geq 1$ are equal and of codimension 
$r$ in
$H^0(K_X^{\otimes \frac{m+2}{2}})$. 

\item{3)} The map $\alpha_n$ surjects if $\frac{m+2}{2} \leq n
\leq m-1$. 

\item{4)} The images of the maps $\alpha(s,t)$, when $s+r=m+1$,
$s,t \geq 1$ are equal and of codimension $r-1$ in
$H^0(K_X^{\otimes m+1})$.

\item{5)} The map $\alpha(\frac{m+2}{2},\frac{m+2}{2})$
 surjects if $r=1$.

 \item{6)} The map $\alpha_{m+1}$ surjects if $r \geq 2$.

  \item{7)} The map $\alpha_n$ surjects if $n \geq m+2$.

   \medskip

     Before we prove the above claims, we will show how they imply the
      result. Claim 1) implies that the part of $R$ of degree less than or
      equal to $\frac{m}{2}$ is generated in degree $1$.
      Claim 2) implies that the subspace of the part of degree
      $\frac{m+2}{2}$ generated in lower degree has codimension $r$, hence,
      in order to generate $R$, $r$ generators of degree $\frac{m+2}{2}$ are
      needed. Claim 1) and 3) imply that the part of $R$ of degree less than
      or equal to $m$ is generated in degree $1$ and $\frac{m+1}{2}$. By the
      same argument as before, Claim 4) implies that, in order to generate
$R$, we need $r-1$ linearly independent elements of degree
$m+1$. Claims 5) and 6) imply that $R_{m+2}$ is generated in degree
$\frac{m+2}{2}$ or lower.  Finally, Claim 7) proves that
the part of $R$ of degree greater than $m+2$ is generated in degree
$m+2$ or lower. Summarizing, the above claims show that $R$ is
generated by its part of degree $1$, by $r$ linearly independent
elements in degree
$\frac{m+2}{2}$ and by $r-1$ linearly independent elements in degree
$m+1$.

Now we proceed to prove claims 1) to 7). The proof goes like the
 proof of \oddmain \ns. Thus Claims 1), 3), 6) and 7) follow from
 \evencurve  (2.1) and
\superlemma (2). The proof of Claim 2) follows from \evencurve
(2.2) and \evencurve (2.5),  and from  \superlemma (3), arguing like for the
proof the 3) in the proof of \oddmain \ns.  Likewise Claim 4) follows from
\evencurve (2.3) and  2.10 (2.6) and from \superlemma \ns (3). 
Finally Claim 5) follows from
\evencurve
(2.3), (2.4)  and 2.10 (2.6) and from \superlemma (3). 
\qed

\proclaim{\evenremark} \evencurve states
 explicitely that the ring 
$\bigoplus_{n \geq
 0}H^0(C,\theta^{\otimes n})$ has in each degree the same number of
 generators as $\bigoplus_{n \geq
 0}H^0(X,K_X^{\otimes n})$. The way \evenmain is proven shows that
 the same is true for the ``intermediate" rings
 $\bigoplus_{n
 \geq 0}H^0(X_{m'},L_{m'}^{\otimes n})$, $2 \leq m' \leq m-1$.
 \endproclaim

  We find now sufficient conditions for $K_X^{\otimes n}$ to satisfy
  property $N_0$: 

  \proclaim{\evenNo}  Let $X$ be a smooth, pluriregular variety of
   general type of even dimension $m$ with ample and
   base-point-free canonical bundle.  Assume that the canonical
   morphism
   has degree $3$ and maps $X$ onto a variety $Y$
   of minimal degree
   $r$. 
   \item{1)} If $r=1$, then the line bundle $K_X^{\otimes n}$ satisfies
     propety $N_0$ if and only if $n \geq \frac{m+2}{2}$.

     \item{2)} If $r > 1$,
      then the line
      bundle
      $K_X^{\otimes n}$ fails to satisfy property
      $N_0$ if $n \leq \frac m 2$ and satisfies property $N_0$
      if $n \geq m+1$.
      \endproclaim

      \noindent {\it Proof.}
       The proof of 1) follows steps similar to those of the proof of \text{\oddNo
       \ns.}

       {\it Step 1.} $K_X^{\otimes n}$ does not satisfy property $N_0$ if $n
        \leq \frac{m}{2}$.

        {\it Step 2.} $K_X^{\otimes n}$ satisfies property $N_0$ if $n =
         \frac{m+2}{2}$.

         {\it Step 3.} $K_X^{\otimes n}$ satisfies  property $N_0$ if
          $\frac{m+4}{2} \leq n
          \leq m+1$.

          {\it Step 4.} $K_X^{\otimes n}$ satisfies  property $N_0$ if
           $n
           \geq m+2$.
Step 1 is true by the same reason as Step 1 of the proof of
\oddNo \ns.
The other steps are proven like in \oddNo and the key facts to use
are the surjectivity of $\alpha(\frac{m+2}{2}, \frac{m+2}{2})$ if $r
=1$ and the surjectivity of $\alpha_{n'}$ for all $n' \geq m+2$. This
was shown in 5) and 7) of the proof of \evenmain \ns. We leave the rest
of the details to the reader.

The proof of 2) goes along the same lines.
 The fact that the line
 bundle
 $K_X^{\otimes n}$ fails to satisfy property
 $N_0$ if $n \leq \frac m 2$ is proven like Step
 1 of \oddNo \ns. The fact that $K_X^{\otimes n}$
 satisfies property $N_0$
 if $n \geq m+1$ follows from the surjectivity of $\alpha_{n'}$ for all
 $n' \geq m+1$ and $r > 1$.
 \qed

   \bigskip
   Comparing the previous result with \oddNo we realize that
   in the even dimensional case the following question is left
open:

\proclaim{\evenNoquestion} Let $X$ be a smooth variety of general
type and dimension $m$ even and let $X @>
\pi >> Y$ be a generically finite, canonically induced, degree
$3$ cover of a variety
$Y$ of minimal degree $r > 1$. Does
$K_X^{\otimes n}$ satisfy property
$N_0$ if $\frac{m+2}{2} \leq n \leq m$?
\endproclaim

The reason why this question cannot be addressed with the present
 arguments
 is that we cannot decide whether the map
 $\alpha(\frac{m+2}{2},\frac{m+2}{2})$ surjects. The map
 $\beta(\frac{m+2}{2}, \frac{m+2}{2})$ does surject, as stated in
 \evencurve (2.4),
 but if we were to use \superlemma to deduce from the surjectivity of
 $\beta(\frac{m+2}{2}, \frac{m+2}{2})$ the
 surjectivity of $\alpha(\frac{m+2}{2},\frac{m+2}{2})$, we would also
 need the surjectivity of $\beta(\frac{m+2}{2}, \frac{m}{2})$. The
 latter map is however non surjective (cf. \evencurve (2.3)).

   \heading
   3. On the structure of the canonical morphism 
\endheading

In the previous section we studied the canonical ring of 
pluriregular varieties of general type $X$ which were triple canonical
covers of varieties of minimal degree.
In this
section we relax our hypothesis on $X$ and 
$\pi$ and study the structure of the
canonical morphism $X @> \pi >> Y$. This study also tells us
about the structure of $X$ and $Y$.  
We do not not assume a priori $X$ to be
pluriregular or even regular, neither do we assume $\pi$ to be
induced by the
complete canonical series of $X$. In fact we start only with a
base-point-free
canonical subseries (we will see later, cf. \flatstruct \ns, that assuming
$\pi$ to be flat  is
sufficient for $X$ to be as the varieties studied in Section 2, that
is, pluriregular
and with $\pi$ induced by the complete canonical series).
Despite these weaker hypothesis,
in the next theorem
we are able to obtain interesting information relating the degree of
$Y$ and the dimension of $X$:

\proclaim{\genstruct}  Let $X$ be a smooth variety of general type
of dimension $m \geq 2$
and let $Y
\subset
\bold P^{r+m-1}$ be an $m$-fold of minimal degree. Assume that
$K_X$ is base-point-free and let $W$ be a sublinear series of $|K_X|$
without base points. Let
$\pi: X \longrightarrow Y$ be the morphism induced by $W$.
Assume that $\pi$ is
generically finite of degree $3$. Then

\item{1)} If $m$ is odd, then
 the degree $r$ of $Y$ is even (in particular $Y$ is not
 linear
 $\bold P^m$)

 \item{2)} If $h^1(\Cal O_X)=0$, then $\pi$ is induced by the complete canonical
  series of $X$.  \endproclaim

  \noindent {\it Proof.} Since $W$ has no base points, we may choose
   $m-1$ general members  of
   $W$ so that its intersection $C$ is a smooth and irreducible curve in
   $X$ and so that $\varphi=\pi|_C$ is finite onto a smooth rational
   normal curve $D$
   of degree $r$ in $\bold P^r$.
   Let $\theta=K_X|_C$. By adjunction
   $\theta^{\otimes m} = K_C$.
   Now, if $m$ is odd, then $C$, $\theta$ and $\varphi$ satisfy the hypothesis
   of \oddcurve \ns. Then according to \oddcurve (1.1), $r$ is even.
   On the other hand, $C$, $\theta$ and $\varphi$ satisfy the hypothesis
   of \oddcurve if $m$ is odd, and the hypothesis of \evencurve
   if $m$ is even. Hence, since $m \geq 2$, according to \oddcurve (1.2)
   and \evencurve (1.2),
    $\varphi$ is
    induced by the complete linear series $|\theta|$. Since $H^1(\Cal O_X)=0$, by
    Kawamata-Viehweg vanishing we have also $H^1(\Cal
    O_{X_{m-1}})=\cdots = H^1(\Cal O_{X_2})=0$. Then the fact that
    $\varphi$ is induced by the complete $|\theta|$ implies that $\pi$
    is induced by the complete linear series of the canonical of $X$. This
    proves 2).
    \qed

    \bigskip

\genstruct tells us among other things that a variety of minimal
degree $Y$ with degree an even integer, cannot occur as the image of a
canonical morphism $X @> \pi >> Y$ of degree $3$ if the dimension of $X$
is odd. In the next result we further eliminate possible targets of
$\pi$. Applied to the setting of Section 2, in which $Y$ was a
variety of minimal degree and $\pi$ has degree $3$, \scroll tells us that
there are no even dimensional canonical, generically triple covers of
smooth rational normal scrolls. \scroll is more general since
it considers generically finite canonical morphisms of arbitrary
odd degree and smooth scrolls fibered not necessarily over $\bold P^1$.

      \proclaim{\scroll} Let $X$ be a smooth 
 variety of general type of
       dimension $m$ with base-point-free canonical bundle $X$. Let $W$ be a
       linear subseries of $|K_X|$ without base points and let $X @>
       \pi >> Y$
       be the morphism induced by $W$. Assume that $\pi$ is generically finite
        of degree $n$ and that $Y$ is a scroll. Then either $n$ is even or
         $m$ is odd. In particular, if $n=3$, $m$ is odd.
         \endproclaim 

         {\it Proof.}
          Let $G$ be a fiber of Y. Let $H$ be
          the class of the
          hyperplane section of Y.
           Let $g$ be a straight line in $G$. 
           Then $g=G \cdot H^{(m-2)}$, i.e, g is
           the complete
           intersection on $Y$ of $G$ and $m-2$ suitable hyperplane sections.
           Likewise,  $\pi^*g$ is the complete intersection on $X$ of the pullback
           by $\pi$ of $m-2$
           suitable hyperplane sections and $\pi^*G$.
           We choose the
           $m-2$ hyperplane sections suitably so that, by Bertini, the intersection of the
           pullback by $\pi$ of $m-2$ of them is a smooth irreducible surface
           $X'$. The image of $X'$ by $\pi$ is a surface $Y$, which is also a
           smooth
            scroll and has dimension
            $2$. Applying adjunction recursively we see that
            $K_{X'}$ is the restriction of $K_X+(m-2)\pi^*(H)$ to $X'$. Then the
             sectional genus of $\pi^*G$  is $t=(K_{X'}+\pi^*g) \cdot \pi^*g$.
           Now we
            know that $K_X=\pi^*H$, so we have that
            $$t=[K_X+(m-2)\pi^*H+\pi^*G]\cdot \pi^*g
            =[(m-1)\pi^*H+\pi^*G] \cdot
            \pi^*g=(m-1)n \ .$$
            Thus, if $m$ is even, $m-1$ is odd, so $n$ has to be even. And if $n$
            is odd,
            $m-1$
            has to be even, hence $m$ has to be odd.
            $\square$

            \bigskip

The two previous results eliminate several possibilities
for the target of the morphism $X @> \pi >> Y$. We end the 
section characterizing, under the stronger hypothesis of $\pi$
being flat, the targets that actually do occur. We see this in the
next theorem, where we give the list of the only possible targets. In
\examples \ns, we show that all the possibilities allowed by
 \flatstruct do exist. \flatstruct also shows that a variety of general
type which admits a flat, triple canonical morphism  onto a variety of
minimal degree is necessarily pluriregular. The
theorem below also describes some features of the morphism $\pi$.

              \proclaim{\flatstruct}
               Let $X$ be a smooth variety of general type of dimension $m \geq
               2$. Assume that
               $K_X$ is ample and base-point-free. Let $W$ be a linear
               subseries 
               of
               $|K_X|$ without base points. Let $\pi$ be the
               morphism induced by $W$. Assume that $\pi$ is flat and has degree
               $3$. Let
               $Y=\pi(X)$ be a variety of minimal degree $r$. Then $Y$ is one of the
               following:

               \medskip

                \item{1)}  $\bold P^m$ with $m \geq 2$ even.

                  \item{2)} A smooth rational normal scroll of odd dimension $m \geq
                  3$ and even degree $r$. In this case $X$ is fibered over $\bold P^1$,
                  with general fiber a smooth variety of general type of dimension
                  $m-1$. The restriction of $\pi$ to the general fiber is finite and
                  flat, induced by the complete canonical series, and maps onto
                  $\bold P^{m-1}$, fitting therefore in case 1) above.

                    \item{3)} A smooth quadric hypersurface of odd dimension $m \geq
                     3$.

                     \medskip

                      In addition, $\pi$ is  in fact induced by the
                       complete canonical linear series of $X$; $\pi_*\Cal O_X=\Cal O_Y
                       \oplus L^{-1}
                       \oplus L^{-2}$ for a line bundle $L$ such that $L^{-2}$ is not
                       effective; $K_Y=L^{-2}(1)$ and $H^1(\Cal
                       O_X)=\cdots=H^{m+1}(\Cal O_X)=0$.
                       \endproclaim

                       {\it Proof.} Since $\pi$ is flat and $X$ is smooth, $Y$ is also smooth.
                        Thus
                        $Y$ must be one of the following:

                        \item{a)} $\bold P^m$;

                         \item{b)} a smooth rational normal scroll;

                          \item{c)} a smooth quadric hypersurface.

                            Since $\pi$ is finite of degree $3$ and flat, $\pi_*\Cal O_X=\Cal O_Y
                             \oplus E$, where $E$ is a locally free sheaf of rank $2$ (cf. \HM \ns).
                             Assume
                             $E$  does not decompose. Then since $\pi$ is induced by a canonical
                             subseries, $K_X = \pi^*(\Cal O_Y(1))$. Hence by projection formula,
                             $\pi_*K_X=\Cal O_Y(1) \oplus E(1)$. On the other hand,  by relative
                             duality 
                             $\pi_*K_X=K_Y
                             \oplus (E^{*}\otimes K_Y)$, hence $K_Y=\Cal O_Y(1)$ which is a
                             contradiction. Then $E=L_1^{-1} \oplus L_2^{-1}$, for line bundles
                             $L_1$ and
                             $L_2$ such that neither $L_1^{-1}$ nor $L_2^{-2}$ has sections,
                             for $X$ is connected. Furthermore
                             $$\pi_*K_X=\Cal O_Y(1)
                             \oplus L_1^{-1}(1) \oplus L_2^{-1}(1)$$
                              on the one hand and  on the other hand,
                              $$\pi_*K_X=K_Y \oplus (K_Y \otimes L_1) \oplus
                              (K_Y \otimes L_2) \ .$$
                              We have already noticed that $K_Y$ cannot be $\Cal O_Y(1)$.
                              Then either $K_Y=L_1^{-1}(1)$ or $K_Y=L_2^{-1}(1)$. We
                              see first what happens if $K_Y=L_1^{-1}(1)$. In this case
                              $K_Y
                              \otimes L_1$ equals either
                              $\Cal O_Y(1)$ or
                              $L_2^{-1}(1)$. The latter leads to a contradiction.
                              On the other hand, if $K_Y=L_1^{-1}(1)$
                              and $K_Y
                              \otimes L_1=\Cal O_Y(1)$, then $L_1 = L_2^{\otimes 2}$ and
                              $K_Y=L_2^{-2}(1)$. The case $K_Y=L_2^{-1}(1)$ is
                              analogous, yielding $L_2 = L_1^{\otimes 2}$ and
                              $K_Y=L_1^{-2}(1)$.
                              Thus we have shown that
                              $$ \flatone \ \ \pi_*\Cal O_X=\Cal O_Y \oplus L^{-1} \oplus
                              L^{-2} \ ,$$
                               for a line bundle $L$ such that $L^{-2}$ is not
                               effective. Moreover, $\pi_*K_X= \Cal O_Y(1) \oplus L^{-1}(1)
                               \oplus L^{-2}(1)$ and

                               $$ \flattwo \ \  K_Y=L^{-2}(1)\ .$$

                                 Now we study further restrictions on 1), 2) and 3). Let
                                  $Y$ be first $\bold P^m$. Since $K_Y=\Cal O_{\bold
                                  P^m}(-m-1)$, by $\flattwo$
                                  $m$ has to be even. Let $Y$ be now a rational normal scroll. Let
                                  us denote by $H$ the divisor class of its hyperplane section and by
                                  $F$ the class of a fiber. Let $r$ be the degree of $Y$. Then
                                  $K_Y=-mH+(r-2)F$. Again by \flattwo $r$ has to be even and $m$ has to
                                  be odd. The only case left is when $Y$ is a smooth quadric
                                  hypersurface of dimension greater than $2$. By Lefschetz Theorem,
                                  Pic$Y=\bold Z$ and is generated by the hyperplane section $Y$. In
                                  addition
                                  $K_Y=\Cal O_Y(-m)$ by adjunction. Hence by \flattwo $m$ is odd.

                                   Since $L^{-2}(1)=K_Y$ and by projection formula
                                   $$H^0(\pi_*\pi^*(\Cal
                                   O_Y(1))=H^0(\Cal O_Y(1)) \oplus H^0(L^{-1}(1)) \oplus H^0(L^{-2}(1)\
                                   ,$$
                                   checking for each case a), b) c) we see that $H^0(\pi_*\pi^*(\Cal
                                   O_Y(1))=H^0(\Cal O_Y(1))$,  hence $\pi$ is in fact induced by the
                                   complete canonical series of $X$. 
                                    To further study  the
                                    structure of $\pi$ in case 2), let us denote by $G$ the inverse
                                    image of a general fiber $F$ of $Y$. We have the following exact
                                    sequence:
                                    $$ 0 \longrightarrow H^0(K_X(-G)) \longrightarrow H^0(K_X)
                                    \longrightarrow H^0(K_X \otimes \Cal O_G) \longrightarrow
                                    H^1(K_X(-G)) \ .$$
                                    Since
                                    $\pi$ is finite, and $K_X(-G)=\pi^*(\Cal O_Y(H-F))$,
                                    $$H^1(K_X(-G))=H^1(\Cal O_Y(H-F)
                                    \oplus L^{-1}(H-F) \oplus
                                    L^{-2}(H-F)) \ .$$ Since $L^{-1}(H-F)=\Cal O_Y(\frac{-m+1}{2}H +
                                    \frac{r-4}{2}F)$ and $L^{-2}(H-F)=\Cal O_Y(-mH +
                                    (r-3)F)$, the first cohomology group of the three line bundle
                                    vanishes.
                                     Hence the complete linear series $|K_X|$ restricts to the complete
                                    linear series $K_G$, which maps onto $F=\bold P^{m-1}$. $G$ is
                                    smooth by Bertini, and since
                                    $K_X$ is ample  and base-point-free,
                                    $G$ must be connected.
                                     The finiteness of the
                                    morphism from $G$ to $F$ and the fact that its degree is $3$
                                    are clear. To see the flatness it suffices to see that
                                    $\pi_*(\Cal O_X)
                                    \otimes \Cal O_G = (\pi|_G)_*(\Cal O_G)$.

                                     Finally, to see the vanishings of cohomology we will deal case by
                                      case. Recall that $L^{-1}=\Cal O_Y(\frac 1 2(K_Y-H))$ and that
                                      $\pi_*(\Cal O_X)=\Cal O_Y
 \oplus  L^{-1}
 \oplus
 L^{-2}$. If $Y=\bold P^m$, then $H^i(\Cal O_X)=0$ for all $1\leq
 i \leq
 m-1$, since intermediate cohomology of line bundles on $Y$
 vanishes. If $Y$ is a smooth rational normal scroll, then
 $$\displaylines{H^i(\Cal O_X)=H^i(\Cal O_Y) \oplus H^i(\Cal O_Y
 (-\frac 1 2 (m+1)H+\frac 1 2 (r-2)F))  \cr \oplus H^i(\Cal
 O_Y(-(m+1)H+(r-2)F))=0}$$ for all
 $1\leq i
 \leq m-1$.  The latter can be seen using Serre duality and pushing
 down to $\bold P^1$. 
 Now, for $Y$ smooth quadric hypersurface of dimension $m$ we have
 the sequence
 $$0 \longrightarrow \Cal O_{\bold P^{m+1}}(-2) \longrightarrow \Cal
 O_{\bold P^{m+1}} \longrightarrow \Cal O_Y \longrightarrow 0 $$ so
 the vanishing of $H^i(\Cal O_Y)$ for all $1 \leq i \leq m-1$
 follows from the vanishing of the intermediate cohomology of line
 bundles on projective space.
 $\square$

\bigskip 

  Finally we show in the following proposition
  that all the possible varieties of minimal degree
  allowed by \flatstruct do actually occur. Note also that, by the
  previous theorem, the varieties of general type constructed in the
  next proposition are pluriregular and satisfy therefore the hypothesis
  of \main \ns. 

  \proclaim{\examples}

   \item{a)} There exist smooth varieties of general type $X$ with
      base-point-free canonical bundle $K_X$ and

      cyclic triple covers $X \longrightarrow \bold P^m$
      induced by
       the canonical morphism of $X$ if
       and only if $m$ is even.

       \item{b)} Let $m \geq 3$ and let $Q \subset \bold P^{m+1}$ be a smooth
          hyperquadric of dimension $m$. There exist smooth varieties of
            general type $X$ with
              base-point-free canonical bundle $K_X$ and
              cyclic triple covers $X \longrightarrow Q$ induced by
               the canonical morphism of $X$ if
               and only if $m$ is odd.

               \item{c)}  Let $m \geq 3$. There exist smooth
                  rational normal scrolls $S$  of dimension $m$ and degree $r$,
                  smooth varieties of
                    general type $X$ with
                      base-point-free canonical bundle $K_X$ and
                      cyclic triple covers $X \longrightarrow S$ induced by
                       the canonical morphism of $X$ if
                       and only if $m$ is odd and $r$ is even. 
                        \endproclaim

                        \noindent {\it Proof.}
                        \flatstruct takes care of the ``only if'' part of a), b) and c). To
                        construct examples of triple covers on $\bold P^m$, $m$ even, take the
                        triple cyclic
                        cover \linebreak 
$X @> \pi >> \bold P^m$ ramified along a smooth
                        divisor of degree $\frac {3(m+2)}{2}$. In such case, $$K_X=\pi^*(K_{\bold P^m}
                        \otimes \Cal O_{\bold P^m}(m+2))=\pi^*\Cal
                        O_{\bold P^m}(1) \ .$$
 Moreover
                       $$H^0(K_X)= H^0(\Cal O_{\bold P^m}(1)) \oplus H^0(\Cal O_{\bold
                          P^m}(-\frac {m}{2})) \oplus H^0(\Cal O_{\bold P^m}(-m-1))=
                          H^0(\Cal O_{\bold P^m}(1)) \ ,$$ hence $\pi$ is induced by the complete
                          canonical series of $X$.
                          Analogously, to construct examples of triple cyclic
                          covers of a hyperquadric $Q$
                          of odd dimension $m \geq 3$ we take the cover $X @> \pi >> Q$ ramified along
                          a smooth divisor in $Q$, complete intersection of a hypersurface of
                          $\bold P^{m+1}$ of degree $\frac {3(m+1)}{2}$ and $Q$.
                          Finally, to construct examples of triple cyclic covers of smooth
                          rational normal scrolls  of odd dimension $m \geq 3$ and even degree
                          $r$, consider a smooth rational scroll
                           $S$ possessing smooth divisors linearly equivalent to $\frac {m+1} {2}
                           H - \frac
                           {r-2}{2} F$, where $H$ is the hyperplane class of $S$ and $F$ is the
                           class of a fiber. Then  we take the cyclic triple cover
                           $X @> \pi >> S$ ramified along
                           a smooth divisor linearly equivalent to $\frac {m+1} {2} H - \frac
                           {r-2}{2} F$.
                             \qed

                            \bigskip

                             {\bf \Greenremark } \main and \examples show that if $m \geq 3$,
                             there does not exist a converse to the theorem of M. Green (cf.
                             \Green \ns) which says that if $X$ is a smooth, pluriregular
                             variety of general type with base-point-free canonical bundle, then
                             the canonical ring of $X$ is generated in degree less than or equal to
                             $m$ if the image
                             $Y$ of the canonical morphism $\pi$ of
                             $X$ is not a variety of minimal degree. Indeed,   the canonical
                             ring of the varieties constructed in \examples  are generated in
                             degree less than or equal to
                             $\frac{m+1}{2}$ if $m$ is odd and $\frac{m+2}{2}$ if $m$ is
                             even. Therefore, $Y$ not being a canonical degree is a sufficient
                             condition for the canonical ring of $X$ to be generated in degree less
                             than or equal to $m$ but it is not a necessary condition.
                             This is in contrast with the situation in dimension $2$, where the
                             converse of the result of Green is true, as it was proved by the
                             authors in \GP \ns.

\heading
References \endheading

\item{\Catone} F. Catanese, {\it Equations of pluriregular varieties of general
type}, Geometry today (Rome, 1984), 47--67, Progr. Math., {\bf 60}, Birkhauser
Boston, 1985. 

\item{\Cattwo} F. Catanese, {\it Commutative algebra methods and equations of
regular surfaces}, Algebraic geometry, (Bucharest, 1982),
68--111, Lecture Notes in Math., {\bf 1056}, Springer Berlin, 1984. 

\item{\Ci} C. Ciliberto, {\it Sul grado dei generatori dell'anello di
    una superficie di tipo generale},
Rend.
Sem. Mat. Univ. Politec. Torino {\bf 41} (1983)

\item{\GP} F.J. Gallego and B.P. Purnaprajna, {\it On the canonical
    ring of covers of surfaces of minimal degree}, preprint AG/0111052.

\item{\Green} M.L. Green, {\it The canonical ring of a variety of
general type}, Duke Math. J. {\bf 49} (1982), 1087--1113.

\item{\HM} D. Hahn and R. Miranda, {\it Quadruple covers of algebraic
varieties}, J. Algebraic Geom. {\bf 8} (1999), 1--30.

\item{\Miranda} R. Miranda, {\it Triple covers in Algebraic Geometry},
  Amer. J. Math. {\bf 107} (1985), 1123--1158.

                             \end